\documentclass[a4j,12pt]{article}
\usepackage{amsthm}
\usepackage{amsmath}
\usepackage{amssymb}
\usepackage{amsfonts}
\usepackage{graphicx}

\newtheorem{lem}{Lemma}[section]

\newtheorem{prp}{Proposition}[section]

\newtheorem{apd}{Lemma A.}

\def\al{{\alpha}}

\def\bga{{\text{\boldmath $\gamma$}}}

\def\bmu{{\text{\boldmath $\mu$}}}

\def\bga{{\text{\boldmath $\gamma$}}}

\def\bmu{{\text{\boldmath $\mu$}}}

\def\De{{\Delta}}

\def\u{{\text{\boldmath $u$}}}
\def\v{{\text{\boldmath $v$}}}

\def\x{{\text{\boldmath $x$}}}

\def\z{{\text{\boldmath $z$}}}

\def\E{{\mathrm E}}

\def\tr{{\rm tr\,}}
\def\diag{{\rm diag\,}}

\def\ah{{\hat a}}

\def\ch{{\hat c}}

\def\xb{{\overline \x}}

\begin{document}
\renewcommand{\theequation}{\thesection.\arabic{equation}}
\renewcommand{\abstractname}{\normalsize \bf Abstract}

\makeatletter
\@addtoreset{equation}{section}
\def\section{\@startsection {section}{1}{\z@}{-3.5ex plus -1ex minus 
-.2ex}{2.3ex plus .2ex}{\bf}}
\def\subsection{\@startsection {subsection}{1}{\z@}{-3.5ex plus -1ex minus 
-.2ex}{2.3ex plus .2ex}{\bf}}
\makeatother

\baselineskip=24pt

\begin{center}
\vspace*{2cm}

{\Large \bf A one-sample location test based on weighted averaging of two test statistics \\
in high-dimensional data}

\vskip12pt
{\large 
MASASHI HYODO$^{1\ast}$ and TAKAHIRO NISHIYAMA$^{2}$
}
\end{center}

\vskip12pt

\begin{quote}
\begin{description}
\baselineskip=14pt
\itemsep=0pt
\item{$^{1}$} Department of Mathematical Information Science, \\
\noindent \hspace{-10pt}Tokyo University of Science,~E-Mail:caicmhy@gmail.com
\item{$^{2}$} School of Business Administration, Senshu University
\end{description}
\end{quote}

\baselineskip=20pt
\begin{quote}
{\it We discuss a one-sample location test that can be used in the case of high-dimensional data. 
For high-dimensional data, the power of Hotelling's test decreases when the dimension is close to the sample size. 
To address this loss of power, some non-exact approaches were proposed, e.g., 
Dempster (1958, 1960), Bai and Saranadasa (1996) and Srivastava and Du (2006). 
In this paper, we focus on Hotelling's test and Dempster's test. 
The comparative merits and demerits of these two tests vary according to the local parameters. 
In particular, we consider the situation where it is difficult to determine which test should be used, 
that is, where the two tests are asymptotically equivalent in terms of local power. 
We propose a new statistic based on the weighted averaging of Hotelling's $T^2$ statistic and Dempster's statistic 
that can be applied in such a situation. 
Our weight is determined on the basis of the maximum local asymptotic power on a restricted parameter space 
that induces local asymptotic equivalence between Hotelling's test and Dempster's test. 
In addition, some good asymptotic properties with respect to the local power are shown. 
Numerical results show that our test is more stable than Hotelling's 
$T^2$ statistic and Dempster's statistic in most parameter settings. 
}
\end{quote}

\baselineskip=24pt

\section{Introduction}

Let $\x_{1},\x_{2},\cdots,\x_{N}$ be $p$ dimensional observation vectors 
from $\mathcal{N}_p(\bmu, \Sigma)$. We consider the following one-sample hypothesis test
\begin{align*}
H_0 : \bmu=\bmu_0 \ \ \mbox{vs.} \ \ H_1 : \bmu \neq \bmu_0.
\end{align*}
To test the hypothesis $H_0$, traditionally Hotelling's test statistic ($T^2$-statistic) is used, which is defined by 
\begin{eqnarray*}
T^2=N(\bar{\x}-\bmu_0)'S^{-1}(\bar{\x}-\bmu_0),
\end{eqnarray*}
where
\begin{eqnarray*}
\bar{\x}=\frac{1}{N}\sum_{i=1}^N\x_i,~S=\frac{1}{n}\sum_{i=1}^N(\x_i-\bar{\x})(\x_i-\bar{\x})',
\end{eqnarray*}
and $n=N-1$. It is well known that under the null hypothesis 
$H_0$, $(N-p)/(np)T^2$ has an $F$-distribution with degrees of freedom $p$ and 
$N-p$. Let the significance level be chosen as $\alpha$ and the threshold be denoted by $F_{p,N-p}(\alpha)$. 
Then Hotelling's test rejects $H_0$ if
\begin{eqnarray*}
\frac{N-p}{np}T^2>F_{p,N-p}(\alpha).
\end{eqnarray*}

However, 
Hotelling's test has the serious defect that the $T^2$ statistic is undefined when the 
dimension of the data is greater than the sample size. 
In subsequent years, a number of improvements on Hotelling's test in the high-dimensional 
setting were discussed, see e.g., Dempster (1958, 1960), Bai and Saranadasa (1998), Srivastava (2007), 
Srivastava and Du (2008). 
In this paper, we focus on Dempster's non-exact test. 
Dempster (1958, 1960) proposed a non-exact test for the hypothesis $H_0$, 
where the dimension $p$ is possibly greater than the sample size $N$. 
Dempster's test statistic (D-statistic) is defined as 
\begin{eqnarray*}
D_n=\frac{(\overline{\x}-\bmu_0)'(\overline{\x}-\bmu_0)}{\tr S}. 
\end{eqnarray*}
However, the exact null distribution of $D_n$ was not derived. 
Therefore, Fujikoshi et al. (2004) proposed an approximate test procedure based 
on the asymptotic normality
\begin{eqnarray}
\sqrt{n}\frac{D_n-1}{\sqrt{2\ah_2/(\ch \ah_1^2)}}
\xrightarrow{d} \mathcal{N}\left(0,1\right),
\end{eqnarray}
under $H_0$, and the assumptions
\begin{eqnarray*}
&{\rm (A1)}&~n,p\to \infty~{\rm with}~\frac{p}{n}\to c\in (0,1),\\
&{\rm (A2)}&~\displaystyle 0< \lim_{p\to\infty} a_i\left(=\lim_{p\to\infty}\frac{\tr\Sigma^i}{p}\right)<\infty,~i=1,\cdots,6.
\end{eqnarray*}
Here, $\ch=p/n$ and 
\begin{eqnarray*}
\widehat{a}_1=\frac{\mbox{tr}S}{p},~\widehat{a}_2=\frac{n^2}{p(n-1)(n-2)}\left(\mbox{tr}S^2-\frac{(\mbox{tr}S)^2}{n}\right)
\end{eqnarray*}
are the unbiased and consistent estimators of $a_1$ and $a_2$. 
Based on the asymptotic normality (1.1), the approximate Dempser's test rejects $H_0$ if
\begin{eqnarray*}
\sqrt{n}\frac{D_n-1}{\sqrt{2\ah_2/(\hat{c} \ah_1^2)}}\ge z(\alpha),
\end{eqnarray*}
where the selected significance level is $\alpha$ and the threshold is denoted by $z(\alpha)$. 

Hotelling's test is powerful when the dimension of the data set is sufficiently small 
as compared with the sample size. 
However, even when $p\leq n$, Hotelling's test is known to perform poorly if $p$ is close to $n$. 
This behavior was demonstrated by Bai and Saranadasa (1996), who studied 
the performance of Hotelling's test under $p,n\to\infty$ with $p/n\to c<1$, and showed that the 
asymptotic power of the test is decreased for large values of $c$. 
In a comparison of the two tests it can be seen that the power of Hotelling's test increases much more slowly than that of 
Dempster's test, as the non-central parameter increases when $c$ is close to one. 
The conclusion drawn from these results is that the comparative merits and demerits of Hotelling's test and Dempster's test 
vary according to the non-central parameter and $c$. 
The contribution of this paper is that a new statistic that possesses both these properties asymptotically is proposed; 
that is, we propose the following statistic which is a weighted average of the $T^2$ statistic and D-statistic:
\begin{eqnarray*}
T(\rho)=\rho\sqrt{n}\left(\frac{T^2}{n}-\frac{p}{n-p}\right)+(1-\rho)\sqrt{n}(D_n-1),
\end{eqnarray*}
where $\rho\in [0,1]$. 
Then, the method used for determining the weight $\rho$ is an important issue. 
In our study, the weight is determined on the basis of the maximum local asymptotic power. 
The only difficulty is that the true optimal weight depends on the true 
mean vector, which is unobservable. 
One method for erasing the information of the true mean vector is to restrict the parameter space 
that induces local asymptotic equivalence between Hotelling's test and Dempster's test. 
This parameter space results in a situation where it is not easy to determine which test may be used. 
Further, the local asymptotic power on this parameter space is evaluated under the condition of a high dimensional framework, that is, 
the sample size and the dimension simultaneously go to infinity under the condition that $p/n\to c\in (0,1)$. 
Large sample asymptotics assume that the dimension $p$ is finite and fixed, while the sample size $N$ grows indefinitely. 
This asymptotic yields a bad approximation in many real-world situations 
where the dimension $p$ is of the same order as the sample size $N$. 
However, it is well known that the high dimensional approximation 
performs well in not only a high dimensional situation, 
but also a large sample situation. 
This fact explains why high dimensional approximation is used. 
We maximize the local asymptotic power and find the optimal weight as a function of $\Sigma$;~then, we 
obtain its consistent estimator. 
We also show that replacing the true optimal weight with a consistent estimator 
makes no difference asymptotically. 
In addition, when the parameter constraint is removed, our statistic is comparable to Hotelling's test and Dempster's test. 
Our test outperforms both tests in terms of local asymptotic power; that is, we can guarantee 
that our test does not have the lowest local asymptotic power among the three tests. 

This paper is organized as follows. 
In Section 2, we introduce the asymptotic property of Hotelling's test and Dempster's test, 
and propose the asymptotically optimal weight $\rho$ for $T(\rho)$ to address the situation where their local asymptotic powers are equal. 
In addition, we give the sufficient condition of a parameter space that allows our test to 
outperform Dempster's test and Hotelling's test in terms of local asymptotic power. 
In Section 3, we investigate the performances of our test through numerical studies. 
The conclusion of our study is summarized in Section 4. 
Some preliminary results and proofs are given in the appendix. 
\section{Description of the weighted averaging test statistic and its asymptotic properties}

In this section, we propose a weighted averaging test statistic of D-statistic and $T^2$-statistic. 
We consider the class of weighted averaging test statistics
\begin{eqnarray*}
\mathcal{T}=\left\{T(\rho)\left|T(\rho)=\rho\sqrt{n}\left(\frac{T^2}{n}-\frac{p}{N-p}\right)\right.+(1-\rho)\sqrt{n}(D_n-1),\rho\in [0,1]\right\}.
\end{eqnarray*}
We note that class $\mathcal{T}$ includes the D-statistic ($\rho=0$) 
and $T^2$-statistic ($\rho=1$). 

First, we propose the optimal weight on the parameter space such that determining the appropriate use of Dempster's 
test and Hotelling's test is difficult, that is, where Dempster's test and Hotelling's test have same local asymptotic power. 
In order to derive the local asymptotic power of a test statistic belonging to class $\mathcal{T}$, 
we assume the conditions (A1), (A2), and
\begin{eqnarray*}
&{\rm (A3)}&~\displaystyle 0< \lim _{n,p \to \infty} n^{1/2}\Delta^{2} <\infty,
0< \lim _{n,p \to \infty}n^{1/2}\Delta_I^{2}<\infty,~0<\lim _{n,p \to \infty} n^{1/2}\Delta_\Sigma^{2} <\infty,
\end{eqnarray*}
where
\begin{eqnarray*}
&&\Delta^{2}=(\bmu-\bmu_0)'\Sigma^{-1}(\bmu-\bmu_0),~\Delta_I^{2}=(\bmu-\bmu_0)'(\bmu-\bmu_0),~\Delta_{\Sigma}^{2}=(\bmu-\bmu_0)'\Sigma(\bmu-\bmu_0). 
\end{eqnarray*}
The following lemma provides the asymptotic normality of $T(\rho)$ under local alternatives. 
\begin{lem}
Assume conditions {\rm (A1)}, {\rm (A2)}, and {\rm (A3)}. 
For any $\rho\in[0,1]$, it holds that
\begin{eqnarray*}
\frac{1}{\sigma(\rho,\ch,\ah_1,\ah_2)}\left[T(\rho)-\sqrt{n}\left\{\rho\frac{\Delta^2}{1-c}+(1-\rho)\frac{\Delta_I^2}{c a_1}\right\}\right]
\xrightarrow{d} \mathcal{N}\left(0,1\right),
\end{eqnarray*}
where 
\begin{eqnarray*}
&&\sigma^2(\rho,c,a_1,a_2)=\rho^2\frac{2c}{(1-c)^{3}}+(1-\rho)^2\frac{2a_2}{ca_1^2}
+2\rho(1-\rho)\frac{2}{1-c}.
\end{eqnarray*} 
\end{lem}
\noindent
{\bf (Proof)}~See, Appendix A.2.

\noindent
Due to Lemma 2.1, the test based on $T(\rho)$ rejects $H_0$ if
\begin{eqnarray}
\frac{T(\rho)}{\sigma(\rho,\ch,\ah_1,\ah_2)}\ge z(\alpha).
\end{eqnarray}
Now, consider the power for testing procedure (2.1). Let
\begin{eqnarray*}
\delta(\rho|\Delta^2,\Delta_I^2,a_1,a_2)=\frac{\rho\Delta^2/(1-c)+\sqrt{n}(1-\rho )\Delta_I^2/(a_1 c)}
{\sqrt{2\rho^2c(1-c)^{-3}+2(1-\rho)^2a_2/(a_1^2 c)+4\rho(1-\rho)(1-c)^{-1}}}.
\end{eqnarray*}
By using asymptotic normality of $T(\rho)$ (Lemma 2.1), we have
\begin{eqnarray}
&&\Pr\left(\frac{T(\rho)}{\sigma(\rho,c,\ah_1,\ah_2)}\ge z(\alpha)\right)\to
\Phi\left(\sqrt{n}\delta(\rho|\Delta^2,\Delta_I^2,a_1,a_2)-z(\alpha)\right)
\end{eqnarray}
under conditions {\rm (A1)}, {\rm (A2)}, and {\rm (A3)}. 

\medskip
\noindent
Our objective is to determine the weight $\rho$ that maximizes the local asymptotic power (2.2). 
Specifically, we assume a restricted parameter space such that 
the local asymptotic power of Hotelling's test 
and of Dempster's test are asymptotically equivalent. 
By using Lemma 2.1, under assumptions (A1)-(A3) and 
\begin{eqnarray*}
(\boldsymbol{\mu},\Sigma)\in\Omega_0&=&\left\{(\boldsymbol{\mu},\Sigma)\left|\frac{\Delta^2}{\Delta_I^2}\right.=\frac{1}{\sqrt{(1-c)a_2}}\right\},
\end{eqnarray*}
it holds that
\begin{eqnarray}
\Pr\left(\frac{N-p}{np}T^2\ge F_{p,N-p}(\alpha)\right)-\Pr\left(\sqrt{n}\frac{D_n-1}{\sqrt{2\ah_2/(\ch \ah_1^2)}}\ge z(\alpha)\right)\to 0;
\end{eqnarray}
this is, their powers are asymptotically equivalent when $(\boldsymbol{\mu},\Sigma)\in \Omega_0$. 
In the following proposition, we obtain the optimal weight on the parameter space $\Omega_0$. 
\begin{prp}
Assume the conditions $(\boldsymbol{\mu},\Sigma)\in \Omega_0$ and {\rm (A1)}-{\rm (A3)}. 
Then, the statistics
\begin{eqnarray}
T(\rho^\ast(c,a_1,a_2))=
\rho^\ast(c,a_1,a_2)\sqrt{n}\left(\frac{T^2}{n}-\frac{p}{N-p}\right)+(1-\rho^\ast(c,a_1,a_2))\sqrt{n}(D_n-1)~
\end{eqnarray}
has maximum local asymptotic power 
\begin{eqnarray*}
\Phi\left(\frac{\sqrt{n(1-c)}\Delta^2}{\sqrt{c(a_1/\sqrt{a_2(1-c)}+1)}}\right)
\end{eqnarray*}
in class $\mathcal{T}$. Here,
\begin{eqnarray*}
\rho^\ast(c,a_1,a_2)=\left(\frac{a_1c \sqrt{a_2 (1-c)}}{a_2 (1-c)^2}+1\right)^{-1}.
\end{eqnarray*}
\end{prp}
\noindent
{\bf (Proof)}~See, Appendix A.3.

\medskip
\noindent
In practice, it is necessary to replace the unknown parameters $a_1$ and $a_2$ in (2.4) with 
their consistent estimators $\ah_1$ and $\ah_2$. 
Nishiyama et al. (2013) provided the following unbiased and 
consistent estimators of $a_1$, $a_2$, $a_3$:
\begin{eqnarray*}
\hat{a}_1&=&\frac{\tr{S}}{p},\\
\hat{a}_2&=&\frac{n^2}{p(n+2)(n-1)} \left \{ \tr{S^2} -\frac{(\tr{S})^2}{n} \right \},\\
\hat{a}_3&=&\frac{n^2}{(n+4)(n+2)(n-1)(n-2)p} \{ n^2\tr{S^3} - 3n\tr{S^2}\tr{S}+ 2(\tr{S})^3 \}.
\end{eqnarray*}
In this study, $\hat{a}_3$ is used (2.8). 
The following lemma shows the asymptotic properties of these estimators. 

\begin{lem}
Assume conditions $\mathrm{(A1)}$ and $\mathrm{(A2)}$. Then, it holds that
\begin{eqnarray*}
\hat{a}_i=a_i+O_p(n^{-1}),~i=1,2,3.
\end{eqnarray*}
\end{lem}
\noindent
{\bf (Proof)}~See, Hyodo et al. (2014). 

\medskip
\noindent
Using Lemma 2.2, we propose an adapted version of (2.4):
\begin{eqnarray}
T(\rho^\ast(\ch,\ah_1,\ah_2))=
\rho^\ast(\ch,\ah_1,\ah_2)\sqrt{n}\left(\frac{T^2}{n}-\frac{p}{N-p}\right)+(1-\rho^\ast(\ch,\ah_1,\ah_2))\sqrt{n}(D_n-1).~
\end{eqnarray}
Further, we denote $\rho^\ast(\ch,\ah_1,\ah_2)$ by $\hat{\rho}^\ast$ and $\rho^\ast(c,a_1,a_2)$ simply by $\rho^\ast$. \\
According to the asymptotic normality of $T(\hat{\rho}^\ast)$ under the null hypothesis $H_0$, we propose the test rejects $H_0$ if
\begin{eqnarray}
\frac{T(\hat{\rho}^\ast)}{\sigma(\hat{\rho}^\ast,\hat{c},\ah_1,\ah_2)}\ge z(\alpha).
\end{eqnarray}
Since $\hat{\rho}^\ast=\rho^\ast+o_p(1)$, 
we obtain the asymptotic power of (2.6) as 
\begin{eqnarray*}
&&\Pr\left(\frac{T(\hat{\rho}^\ast)}{\sigma(\hat{\rho}^\ast,c,\ah_1,\ah_2)}\ge z(\alpha)\right)
\to \Phi\left(\frac{\sqrt{n(1-c)}\Delta^2}{\sqrt{c(a_1/\sqrt{a_2(1-c)}+1)}}\right).
\end{eqnarray*}
Thus, the power of $T(\hat{\rho}^\ast)$ is asymptotically equivalent to that of 
$T(\rho^\ast)$. 
\medskip

From Proposition 2.1 and the above results, 
we derive the asymptotic null distribution of the proposed test statistic $T(\hat{\rho}^\ast)$; 
the improved estimator of the critical point of our test is derived by using the Cornish-Fisher expansion. 
The following proposition provides the asymptotic null distribution of $T(\hat{\rho}^\ast)/\sigma(\hat{\rho}^\ast,\ch,\ah_1,\ah_2)$. 
\begin{prp}
Assume assumptions {\rm (A1)} and {\rm (A2)} and $H_0$. Then, it holds that
\begin{eqnarray*}
\Pr\left(\frac{T(\hat{\rho}^\ast)}{\sigma(\hat{\rho}^\ast,\ch,\ah_1,\ah_2)}\leq x\right)&=&\Phi(x)-\frac{\phi(x)}{\sqrt{n}}
\left\{\frac{b_1(c)(it)}{\sigma(\rho^\ast,c,a_1,a_2)}\right.\\
& &\left.+\frac{b_3(\rho^\ast,c,a_1,a_2,a_3)(it)^3}{\sigma^3(\rho^\ast,c,a_1,a_2)}(x^2-1)\right\}+o\left(n^{-1/2}\right),
\end{eqnarray*}
where
\begin{eqnarray*}
b_1(c)&=&\nu_1(c),~b_3(c,a_1,a_2,a_3)
=\frac{\nu_3(c,a_1,a_2,a_3)}{6}
-\frac{\nu_1(c)}{2}.
\end{eqnarray*}
Here,
\begin{eqnarray*}
\nu_1(c)&=&\frac{2 \rho^\ast c}{(1-c)^2},\\
\nu_3(c,a_1,a_2,a_3)&=&
\frac{4 \rho^{\ast^3} c(5 c+2)}{(1-c)^5}
+\frac{24 \rho^{\ast^2} (1-\rho^\ast) (c+1)}{(1-c)^3}
+\frac{12 \rho^\ast(1-\rho^\ast)^2 a_2(2-c)}{a_1^2 (1-c)^2 c}\\
& &+\frac{8 (1-\rho^\ast)^3 a_3 }{a_1^3 c^2}.
\end{eqnarray*}
\end{prp}
\noindent
{\bf (Proof)}~See, Appendix A.4. 

\medskip
\noindent
Let $x(\al)$ be the upper $100\al$ percentile of the statistic $T(\hat{\rho}^\ast)/\sigma(\hat{\rho}^\ast,\ch,\ah_1,\ah_2)$. 
In addition, the Cornish-Fisher expansion of the true upper $100\al$ percentile 
is obtained by 
\begin{eqnarray}
x(\al)\approx z(\al)
+\frac{1}{\sqrt{n}}\left\{\frac{a_1(c)}{\sigma(\rho^\ast,c,a_1,a_2)}
+\frac{a_3(\rho^\ast,c,a_1,a_2,a_3)}{\sigma^3(\rho^\ast,c,a_1,a_2)}(z(\al)^2-1)\right\}.
\end{eqnarray} 
In practice, it is necessary to replace the unknown parameters $a_1$, $a_2$, and $a_3$ in (2.7) with 
their consistent estimators $\ah_1$, $\ah_2$, and $\ah_3$. 
We replace the $a_i$'s in (2.7) with their unbiased and consistent estimator $\hat{a}_i$, 
and propose an approximate upper 100$\al$-percentile
\begin{eqnarray}
\widehat{x}(\al)=z(\al)
+\frac{1}{\sqrt{n}}\left\{\frac{a_1(\ch)}{\sigma(\hat{\rho}^\ast,\ch,\ah_1,\ah_2)}
+\frac{a_3(\hat{\rho}^\ast,\ch,\ah_1,\ah_2,\ah_3)}{\sigma^3(\hat{\rho}^\ast,\ch,\ah_1,\ah_2)}(z(\al)^2-1)\right\}.
\end{eqnarray} 
Applying (2.8), the test rejects $H_0$ if
\begin{eqnarray}
\frac{T(\hat{\rho}^\ast)}{\sigma(\hat{\rho}^\ast,\ch,\ah_1,\ah_2)}\ge \widehat{x}(\al).
\end{eqnarray}
\medskip

Finally, we compare Hotelling's test and Dempster's test with our test {\rm (2.6)}(or {\rm (2.9)}). 
In the following proposition, we give the sufficient condition that allows our test to have the highest 
local asymptotic power among the three tests. 
Furthermore, even when a sufficient condition does not hold, we can guarantee 
that our test does not have the lowest local asymptotic power among the three tests. 
\begin{prp}
Assume {\rm (A1)},{\rm (A2)}, and {\rm (A3)}. 
The proposed test {\rm (2.6)}{\rm (}or {\rm (2.9)}{\rm )} has the highest local asymptotic power 
among the three tests under the condition
\begin{eqnarray*}
{\rm (C1)}~\frac{\Delta_{\Sigma^{-1}}^2}{\Delta_I^2}\in
\left[\frac{\sqrt{2}\left(1+a_1\sqrt{(1-c)/a_2}\right)^{1/2}-1}{\sqrt{a_2(1-c)}},
\frac{\left\{\sqrt{2}\left(1+a_1\sqrt{(1-c)/a_2}\right)^{1/2}-1\right\}^{-1}}{\sqrt{a_2(1-c)}}\right].
\end{eqnarray*}
Furthermore, the local asymptotic power of our test {\rm (2.6)}{\rm (}or {\rm (2.9)}{\rm )} is second highest among the three tests 
when condition {\rm (C1)} does not hold. 
\end{prp}
\noindent
{\bf (Proof)}~See, Appendix A.5. 
\section{Numerical results}
In this section, we investigate the finite sample behavior of the proposed test and 
compare it with the $T^2$ test and Dempster's test. 
To compare the three tests, we need to define the Attained Significance Level (ASL) and the empirical
powers. We draw an independent sample of size $N= 40 i+ p$, 
where $i=1,\dots,10$ valid $p$-dimensional normal distributions $\mathcal{N}_p(\bmu,\Sigma)$ 
under the null hypothesis $H_0:\bmu=\boldsymbol{0}$. 
Further, we set the covariance structures $\Sigma=(\eta^{|i-j|})$, where $\eta=0.2,0.4,0.6$, 
respectively. 
We replicate this $r=10^5$ times, and, using $T^2$, 
$D_n$, and $T(\hat{\rho}^\ast)$, calculate 
\begin{eqnarray*}
ASL_{\alpha}\left ( T^2 \right )
= \frac{\sharp ~\textrm{of} \, \left((N-p)/(np)T^2>F_{p,N-p}({\alpha}) \right  )}{r},\\
ASL_{\alpha}\left ( D_n \right )
= \frac{\sharp ~\textrm{of} \, \left(\sqrt{n}(D_n-1)/\sqrt{2\hat{a}_2/(\hat{c}\hat{a}_1^2)} > y({\alpha})  \right  )}{r},
\end{eqnarray*}
and
\begin{eqnarray*}
ASL_{\alpha} \left (T(\hat{\rho}^\ast) \right )=
\frac{\sharp ~\textrm{of} \, \left(T(\hat{\rho}^\ast)/\sigma(\hat{\rho}^\ast,\hat{c},\hat{a}_1,\hat{a}_2) > \widehat{x}({\alpha}) \right  )}{r},
\end{eqnarray*}
denoting the ASL of $T^2$, $D_n$, and $T(\hat{\rho}^\ast)$, respectively. 
Here, $y(\alpha)$ is the improved estimator of the critical point for $D_n$, which 
was provided by Nishiyama et al. (2013) and defined as
\begin{eqnarray*}
y(\alpha) = z(\alpha) + {1 \over \sqrt{p}}q_1(z(\alpha)) 
+ {1 \over p}q_2(z(\alpha)) + {1 \over n}q_3(z(\alpha)),
\end{eqnarray*}
and 
\begin{eqnarray*}
q_1(z(\alpha)) &=& {\sqrt{2}\hat{a}_3 \over 3 \sqrt{\hat{a}_2^3}}(z({\alpha})^2 - 1), \\
q_2(z(\alpha)) &=& {\hat{a}_4 \over 2\hat{a}_2^2}z(\alpha)(z(\alpha)^2-3)
- {2\hat{a}_3^2 \over 9\hat{a}_2^3}z(\alpha)(2z(\alpha)^2-5), \\
q_3(z(\alpha)) &=& {z(\alpha) \over 2},
\end{eqnarray*}
where $\hat{a}_4$ is the consistent estimator of $a_4$. For details, see Nishiyama et al. (2013). 

The attained significance levels specified by the selection of set $(p,\eta)$ are given in Tables 1-6. 
Since Hotelling's test is an exact test under the multivariate normality assumptions, we focus on Dempster's test and our test. 
Tables 1-6 show that the attained significance levels of both tests 
approximate the nominal level $\alpha$ reasonably well in all cases. 
In addition, we note that, according to these results, our test has a tendency to become conservative. 
To compute the empirical powers, we select
\begin{eqnarray*}
\bmu=\left(\frac{2}{n^{1/4}\sqrt{p}},\ldots,\frac{2}{n^{1/4}\sqrt{p}}\right).
\end{eqnarray*}
Using the same number of replications as above, 
we draw independent samples of size $N$ from $\mathcal{N}_p(\bmu,\Sigma)$, and calculate the empirical power as 
\begin{eqnarray*}
EP_{\alpha}\left ( T^2 \right )
= \frac{\sharp ~\textrm{of} \, \left((N-p)/(np)T^2>F_{p,N-p}({\alpha}) \right  )}{r},\\
EP_{\alpha}\left ( D_n \right )
= \frac{\sharp ~\textrm{of} \, \left(\sqrt{n}(D_n-1)/\sqrt{2\hat{a}_2/(\hat{c}\hat{a}_1^2)} > y({\alpha})  \right  )}{r},
\end{eqnarray*}
and
\begin{eqnarray*}
EP_{\alpha}\left (T(\hat{\rho}^\ast) \right )=
\frac{\sharp ~\textrm{of} \, \left(T(\hat{\rho}^\ast)/\sigma(\hat{\rho}^\ast,\hat{c},\hat{a}_1,\hat{a}_2) > \widehat{x}({\alpha}) \right  )}{r}.
\end{eqnarray*}
The results for the empirical power are summarized in Tables 7 to 12, 
where bold face marks the highest power among the three tests. 
These tables show that, while our test statistic has the highest power among the three tests in many cases, the other 
two tests have the highest power in some cases. 
Specifically, among the three tests, the performance of Dempster's test is 
comparatively good when $N$ is small, and that of Hotelling's test is comparatively good when $N$ is large. 
Although the power of our test is not always the highest, it is close to being so. 
In other words, the weight behaves such that our statistic is comparable with whichever statistic has 
the relatively higher power, the D-statistic or the $T^2$-statistic. 
\section{Conclusion}
We proposed a new test statistic for the one-sample location test in high-dimensional data. 
Our proposed test statistic uses the weighted averaging of Hotelling's $T^2$ statistic and Dempster's statistic. 
Some asymptotic properties of this statistic were also shown. 
The important issue is that the local asymptotic power of our test does not become lower than that 
of Hotelling's test and Dempster's test. 
In addition, simulations indicate that the newly derived test statistic is relatively stable as compared with the 
D-statistic and $T^2$-statistic. 
When the difference in the power of the D-statistic and the $T^2$-statistic is large, 
it can be seen that our statistic is comparable with whichever statistic has the 
relatively higher power, the D-statistic or $T^2$-statistic. 
In conclusion, we recommend that our test statistic be applied instead of the D-statistic and $T^2$-statistic over a wide range. 
\medskip

\noindent
\section*{Acknowledgments}
The authors thank Professor Yasunori Fujikoshi and Professor Takashi Seo 
for extensive discussions, references, and encouragement, 
and Ms. Manami Okuyama for the Monte Carlo simulation used to obtain the numerical results. 
\newpage
%
%
\section*{Appendix A.}

 \def\theequation{A.\arabic{equation}}
  \makeatletter
  \makeatother

\subsection*{A.1.~~Some preliminary result}
\begin{apd}[The central limit theorem for quadratic forms]
Let $\z=(z_1,\cdots, \\ z_p)'$ be distributed $p$-dimensional standerd normal random variable and $\Omega=
{\rm diag}(\omega_1,\cdots,\\ \omega_p)$ 
be arbitrary $p\times p$ non random diagonal matrix. 
Suppose that $T=\z'\Omega\z-\tr\Omega$ and 
$\sigma_p^2=2\tr\Omega^2$. 
Then, $T/\sigma_p\xrightarrow{d}\mathcal{N}(0,1)$ as $p\to\infty$ if the following condition is satisfied:
\begin{eqnarray}
\frac{\tr\Omega^4}{(\tr\Omega^2)^2}\to 0~{as}~p\to\infty.
\label{con}
\end{eqnarray}
\end{apd}
\noindent
{\bf (Proof)}

\noindent
It can be expressed that
\begin{eqnarray*}
T&=&(\z'\Omega\z-\tr\Omega)\\
&=&\sum_{i=1}^n(\omega_{i}z_{i}^2-\omega_{i}).
\end{eqnarray*}
Let $Y_i=\omega_{i}z_{i}^2-\omega_{i},~i=1,2,\ldots,p$. 
Then $T=\sum_{i=1}^nY_i$ and the moment of $Y_i$ is caluclated by
\begin{eqnarray*}
\E[Y_i^2]&=&2 \omega_{i}^2,~\E[Y_i^4]=60\omega_{i}^4.
\end{eqnarray*}
We wish to give sufficient conditions that ensure $T/\sigma_p\xrightarrow{d}\mathcal{N}(0,1)$. 
For now, we check only the Lyapunov Condition. 
The the Lyapunov Condition for sequences $\{Y_i\}_{i=1}^p$ states that
\begin{eqnarray*}
{\rm there~exists}~\eta\in\mathbb{N}~{\rm such~that}~\frac{\sum_{i=1}^p\E[Y_i^{2+\eta}]}{\sigma_p^{2+\eta}}\to 0~{\rm as}~p\to\infty.
\end{eqnarray*}
Based on the first and second moments of $Y_i$, we can caluclate
\begin{eqnarray}
\sum_{i=1}^p\E[Y_i^2]=2\tr\Omega^2(\equiv \sigma_p^2),~\sum_{i=1}^p\E[Y_i^4]=60\tr\Omega^4\label{eq2}.
\end{eqnarray}
From (\ref{eq2}), under the condition (\ref{con}),
\begin{eqnarray*}
\frac{\sum_{i=1}^p\E[Y_i^4]}{\sigma_p^4}
=\frac{60\tr\Omega^4}{(2\tr\Omega^2)^2}\to 0
\end{eqnarray*}
as $p\to\infty$. This result show that the condition (\ref{con}) implie Lyapunov Condition. 
Thus, the Lyapunov Condition also implies $T/\sigma_p\xrightarrow{d}\mathcal{N}(0,1)$. $\hfill\square$
\begin{apd}[Some moments for quadratic forms]
Let $\z$ be distributed $p$-dimensi\\onal standerd normal random variable and $A_i,~i=1,2,3$ 
be arbitrary $p\times p$ diagonal matrix. Then it holds that
\begin{eqnarray*}
&&{(\rm i)}\E[\z'A_1\z]=\tr A_1,\\
&&{(\rm ii)}\E[\z' A_1\z\z' A_2\z]=2\tr A_1A_2+\tr A_1\tr A_2,\\
&&{(\rm iii)}\E[\z' A_1\z\z' A_2\z\z'A_3\z]=
\tr A_1 \tr A_2 \tr A_3 + 2 \tr A_3 \tr A_1A_2 \\
&&~~~~~~~~~~~~~~~~~~~~~~~~~~~~~~~~~~~~
+2 \tr A_2\tr A_1A_3 + 2 \tr A_1\tr A_2A_3+ 
8 \tr A_1A_2A_3.
\label{con2}
\end{eqnarray*}
\end{apd}
\noindent
{\bf (Proof)}~See e.g. Mathai et al. (1995). 
\subsection*{A.2.~~Proof of Lemma 2.1.}
At first, we expand  $T^2$ stochastically. 
Suppose that $\Gamma=(\bga_1,\cdots,\bga_p)$ is an orthogonal matrix such that 
$\Sigma=\Gamma\Lambda\Gamma'$, where $\Lambda=\diag(\lambda_1,\ldots,\lambda_p)$ and, for 
$i=1,\ldots,p$, $\lambda_i$ is $i$-th eigenvalue of $\Sigma$. 
Define the random variables $\u$ and $W$ by
\begin{eqnarray*}
\u=&
\sqrt{N}\Gamma'(\Sigma^{-1/2}(\xb-\bmu_0)-\boldsymbol{\tau})
,~W=n\Gamma'\Sigma^{-1/2}S\Sigma^{-1/2}\Gamma.
\end{eqnarray*}
It is seen that $\u$ and $W$ are mutually independently distributed as 
$\u\sim\mathcal{N}_p(\boldsymbol{0},I_p)$, respectively, 
where $\boldsymbol{\tau}=\Sigma^{-1/2}(\boldsymbol{\mu}-\boldsymbol{\mu}_0)$. 
Then the statistic $T^2/n$ is denoted by 
\begin{equation*}
{(\Gamma\u+\sqrt{N}\boldsymbol{\tau})}^{'}W^{-1}{(\Gamma\u+\sqrt{N}\boldsymbol{\tau})}\stackrel{d}=
\frac{(\Gamma\u+\sqrt{N}\boldsymbol{\tau})'(\Gamma\u+\sqrt{N}\boldsymbol{\tau})}
{\boldsymbol{v}'\boldsymbol{v}},
\end{equation*}
where $\boldsymbol{v}\sim\mathcal{N}_{N-p}(\boldsymbol{0},I_{N-p})$, 
and $\boldsymbol{u}$ and $\boldsymbol{v}$ are mutually independent. 
Then the the statistic $T^2/n$ can be expanded as 
\begin{eqnarray*}
\frac{T^2}{\sqrt{n}}
&\stackrel{d}=&
\sqrt{n}\frac{(\boldsymbol{u}'\boldsymbol{u}-p)+2\sqrt{N}\boldsymbol{\tau}'\Gamma\boldsymbol{u}+N\boldsymbol{\tau}'\boldsymbol{\tau}+p}
{(N-p)}\left(1+\frac{\boldsymbol{v}'\boldsymbol{v}-(N-p)}{N-p}\right)^{-1}\\
&=&\sqrt{n}\left(\frac{N\Delta^2}{N-p}+\frac{p}{N-p}\right)+
\frac{\boldsymbol{u}'\boldsymbol{u}-p}{\sqrt{n}(1-c)}
-\frac{c}{1-c}\frac{\boldsymbol{v}'\boldsymbol{v}-n(1-c)}{\sqrt{n}(1-c)}+o_p(1).
\end{eqnarray*}
Thus, we have 
\begin{eqnarray}
\sqrt{n}\left\{\frac{T^2}{n}-\left(\frac{N\Delta^2}{N-p}+\frac{p}{N-p}\right)\right\}
&\stackrel{d}=&
\frac{\boldsymbol{u}'\boldsymbol{u}-p}{\sqrt{n}(1-c)}
-\frac{c}{(1-c)^2}\frac{\boldsymbol{v}'\boldsymbol{v}-n(1-c)}{\sqrt{n}}+o_p(1).\nonumber\\
\end{eqnarray}
\medskip

Next, we expand $D_n$ stochastically as following
\begin{eqnarray*}
\frac{N(\boldsymbol{\bar{x}}-\boldsymbol{\mu}_0)'
(\boldsymbol{\bar{x}}-\boldsymbol{\mu}_0)}{\tr S}
&\stackrel{d}=&\frac{\u'\Lambda\u
+N\boldsymbol{\tau}'\Sigma\boldsymbol{\tau}}{\tr\Sigma}
+o_p(n^{-1/2})\\
&=&\left(1+\frac{\Delta_I^2}{c a_1}\right)
+\frac{\u'\Lambda\u-pa_1}{pa_1}
+o_p(n^{-1/2}).
\end{eqnarray*}
Thus, we can obtain
\begin{eqnarray}
\sqrt{n}\left\{D-\left(1+\frac{\Delta_I^2}{c a_1}\right)\right\}
&\stackrel{d}=&
\frac{\u'\Lambda\u-pa_1}{\sqrt{n}ca_1}
+o_p(1).
\end{eqnarray}
From (A.3) and (A.4), we expand $T(\hat{\rho}^\ast)$ stochastically as following
\begin{eqnarray*}
T(\rho)&=&\rho\left(
\frac{\boldsymbol{u}'\boldsymbol{u}-p}{\sqrt{n}(1-c)}
-\frac{c}{(1-c)^2}\frac{\boldsymbol{v}'\boldsymbol{v}-n(1-c)}{\sqrt{n}}\right)
+(1-\rho)
\frac{\u'\Lambda\u-pa_1}{\sqrt{n}ca_1}+o_p(1)\\
&=&\boldsymbol{u}'\left(\frac{\rho}{\sqrt{n}(1-c)} I_p+\frac{1-\rho}{\sqrt{n}ca_1}\Lambda\right)\boldsymbol{u}
-\tr\left(\frac{\rho}{\sqrt{n}(1-c)} I_p+\frac{1-\rho}{\sqrt{n}ca_1}\Lambda\right)\\
& &-\left\{\boldsymbol{v}'\left(\frac{c\rho}{\sqrt{n}(1-c)^2}I_{N-p}\right)\boldsymbol{v}
-\tr\left(\frac{c\rho}{\sqrt{n}(1-c)^2}I_{N-p}\right)\right\}+o_p(1).
\end{eqnarray*}
By using Lemma A.1 and the independency $\u$ and $\v$, we obtain Lemma 2.1. $\hfill\square$
\subsection*{A.3.~~Proof of Proposition 2.1.}
Assume that $(\boldsymbol{\mu},\Sigma)\in \Omega_0$ i.e. $\Delta_I^2/\Delta^2=\sqrt{(1-c)a_2}$. 
By using Lemma 2.1, we have
\begin{eqnarray*}
&&\Pr\left(\frac{T(\rho)}{\sigma(\rho,\ch,\ah_1,\ah_2)}\ge z(\alpha)\right)
\to \Phi\left(\sqrt{n}f(\rho)\Delta^2-z(\alpha)\right),
\end{eqnarray*}
where
\begin{eqnarray*}
f(\rho)=\frac{\rho/(1-c)+(1-\rho)\sqrt{a_2(1-c)}/(a_1 c)}
{\sqrt{2\rho^2c/(1-c)^3+2(1-\rho)^2a_2/(a_1^2c)
+4\rho(1-\rho)/(1-c)}}.
\end{eqnarray*}
To obtain the optimal $T(\rho)$ which maximize the local asymptotic 
power function, we consider the optimization problem:
$\max_{\rho\in [0,1]} f(\rho)$, 
because $\Phi(\cdot)$ is monotonically increasing on $\mathbb{R}$. 
We find $f'(\rho)$ and set it equal to zero. 
Solving $f'(\rho)=0$ for $\rho$ gives us
\begin{eqnarray*}
\rho^\ast(c,a_1,a_2)=\left(\frac{a_1c \sqrt{a_2 (1-c)}}{a_2 (1-c)^2}+1\right)^{-1}.
\end{eqnarray*}
The second derivative is given by
\begin{eqnarray*}
& &\hspace{-15pt}f''(\rho)\\
&=&
-\frac{\Delta^2
\left\{
(1-c)^2\left(\sqrt{a_2(1-c)}-a_2/a_1\right)
+a_1 c (1-c)-c\sqrt{a_2(1-c)}\right\}^4
}
{
\sqrt{2}a_1 (1-c)^3c
\left\{c\ell
\left(a_1^2(1-c)\ell
+\left(\sqrt{a_2(1-c)}-a_1(1-c)\right)^2\right)
\right\}^{3/2}}<0,
\end{eqnarray*}
so $f(\rho^\ast(c,a_1,a_2))$ is a local maximum value. Here, $\ell=a_2/a_1^2+c-1$. 
Since $f'(\lambda)$ is monotone decreasing function on $[0,1]$, 
we can get $\rho^\ast(c,a_1,a_2)$ as the solution to $\max_{\rho\in [0,1]} f(\rho)$. 
Thus, the optimal linear combination is given by
\begin{eqnarray*}
T(\rho^\ast(c,a_1,a_2))=
\rho^\ast(c,a_1,a_2)\sqrt{n}\left(\frac{T^2}{n}-\frac{p}{N-p}\right)+(1-\rho^\ast(c,a_1,a_2))\sqrt{n}(D_n-1)
\end{eqnarray*}
and its asymptotic power is
\begin{eqnarray*}
&&\Pr\left(\frac{T(\rho^\ast(c,a_1,a_2))}{\sigma(\rho^\ast(c,a_1,a_2),c,\ah_1,\ah_2)}\ge z(\alpha)\right)
\to \Phi\left(\frac{\sqrt{n(1-c)}\Delta^2}{\sqrt{c(a_1/\sqrt{a_2}(1-c)^{1/2}+1)}}-z(\alpha)\right).\hfill\square
\end{eqnarray*}
\subsection*{A.4.~~Proof of Proposition 2.2.}
Define the variables 
\begin{eqnarray*}
h_1&=&\frac{\u'\u-p}{\sqrt{2p}},
h_2=\frac{\v'\v-(N-p)}{\sqrt{2(N-p)}},
h_3=\frac{\u'\Lambda\u-pa_1}{\sqrt{2a_2p}},
h_4=\frac{\sqrt{np}(\ah_1-a_1)}{\sqrt{2a_2}}.
\end{eqnarray*}
Since $\hat{a}_i=a_i+o_p(n^{-1/2})$ for $i=1,2$, it holds that $\hat{\rho}^\ast=\rho^\ast+o_p(n^{-1/2})$. 
Thus we obtain
\begin{eqnarray*}
T(\hat{\rho}^\ast)=
\rho^\ast\left(T_1+\frac{T_2}{\sqrt{n}}\right)
+(1-\rho^\ast)\left(D_1+\frac{D_2}{\sqrt{n}}\right)
+o_p(n^{-1/2}),
\end{eqnarray*}
where
\begin{eqnarray*}
T_1&=&\frac{c}{1-c}\left(\frac{\sqrt{2}h_1}{\sqrt{c}}-\frac{\sqrt{2}h_2}{\sqrt{1-c}}\right),
T_2=\frac{c}{1-c}\left(\frac{2h_2^2}{1-c}-\frac{2h_1h_2}{\sqrt{1-c}\sqrt{c}}\right),\\
D_1&=&\frac{\sqrt{2a_2}h_3}{a_1\sqrt{c}},
D_2=-\frac{\sqrt{2a_2}h_4}{a_1\sqrt{c}}.
\end{eqnarray*}
Then the moment of order $i(=1,2,3)$ of $T(\rho^\ast)$ denotes 
\begin{eqnarray*}
\E[T(\hat{\rho}^\ast)]&=&
\E\left[\frac{\sqrt{2c}\rho^\ast h_1}{1-c}
-\frac{\sqrt{2}c\rho^\ast h_2}{(1-c)^{3/2}}
+\frac{\sqrt{2a_2} (1-\rho^\ast)h_3}{a_1\sqrt{c}}\right.\\
& &+\frac{1}{\sqrt{n}}\E\left[
\frac{2c\rho^\ast h_2^2}{(1-c)^2}-\frac{2\rho^\ast\sqrt{c}h_1h_2}{(1-c)^{3/2}}
-\frac{\sqrt{2a_2}(1-\rho^\ast)h_4}{a_1\sqrt{c}}
\right]+o(n^{-1/2}),\\
\E[T^2(\hat{\rho}^\ast)]&=&
\E\left[
\frac{2 \rho^{\ast^2}ch_1^2}{(1-c)^2}
+\frac{2 \rho^{\ast^2} c^2 h_2^2}{(1-c)^3}
+\frac{2(1-\rho^{\ast})^2a_2h_3^2}{a_1^2 c}
-\frac{4 \rho^{\ast^2}c^{3/2}h_1 h_2}{(1-c)^{5/2}}\right.\\
& &+\left.\frac{4\rho^{\ast}(1-\rho^{\ast})\sqrt{a_2}h_1 h_3}{a_1 (1-c)}
-\frac{4 \rho^{\ast}(1-\rho^{\ast})\sqrt{c a_2}h_2 h_3}{a_1 (1-c)^{3/2}}\right]
+o(n^{-1/2}),\\
\E[T^3(\hat{\rho}^\ast)]&=&
\E\left[\frac{2 \sqrt{2} \rho^{\ast^3} c^{3/2} h_1^3}{(1-c)^3}
-\frac{2 \sqrt{2} \rho^{\ast^3} c^3 h_2^3}{(1-c)^{9/2}}
+\frac{6 \sqrt{2c a_2} \rho^{\ast^2}(1-\rho^{\ast})h_1^2h_3}{a_1(1-c)^2}\right.\\
& &
+\frac{6 \sqrt{2}a_2 \rho^{\ast}(1-\rho^{\ast})^2h_1 h_3^2}{a_1^2 (1-c) \sqrt{c}}
+\frac{2 \sqrt{2}a_2^{3/2}(1-\rho^{\ast})^3h_3^3}{a_1^3c^{3/2}}
+\frac{1}{\sqrt{n}}
\left(\frac{12 \rho^{\ast^3}c^3 h_2^4}{(1-c)^5}\right.\\
& &
+\frac{36\rho^{\ast^3}c^2 h_1^2 h_2^2}{(1-c)^4}
-\frac{6 \sqrt{2 a_2} c^{3/2} \rho^{\ast^2}(1-\rho^{\ast})h_2^2h_4}{a_1 (1-c)^3}
+\frac{48\sqrt{a_2}c \rho^{\ast^2}(1-\rho^{\ast})h_1 h_2^2 h_3}{a_1 (1-c)^3}\\
& &
+\left.\left.\frac{12a_2\rho^{\ast}(1-\rho^{\ast})^2h_2^2h_3^2}{a_1^2 (1-c)^2}\right)\right]
+o(n^{-1/2}).
\end{eqnarray*}
By using Lemma A.2, we have
\begin{eqnarray*}
&&\E[h_1]=0,\E[h_2]=0,\E[h_3]=0,\E[h_4]=0,
\E[h_1^2]=1,\E[h_2^2]=1,\E[h_3^2]=1,\E[h_4^2]=1,\\
&&\E[h_1^3]=\frac{2\sqrt{2}}{\sqrt{p}},\E[h_2^3]=\frac{2\sqrt{2}}{\sqrt{N-p}},
\E[h_3^3]=\frac{2\sqrt{2}\tr\Sigma^3}{(\tr\Sigma^2)^{3/2}},
\E[h_1^2h_3]=\frac{2\sqrt{2}\tr\Sigma}{p\sqrt{\tr\Sigma^2}},\\
&&\E[h_1h_3^2]=\frac{2\sqrt{2}}{\sqrt{p}}.
\end{eqnarray*}
Hence, the moments can be calculated by
\begin{eqnarray}
\E[T(\hat{\rho}^\ast)]&=&\frac{\nu_1(c)}{\sqrt{n}}+o\left(n^{-1/2}\right),\\
\E[T^2(\hat{\rho}^\ast)]&=&\sigma^2(c,a_1,a_2)+o\left(n^{-1/2}\right),\\
\E[T^3(\hat{\rho}^\ast)]&=&\frac{\nu_3(c,a_1,a_2,a_3)}{\sqrt{n}}+o\left(n^{-1/2}\right).
\end{eqnarray}
The relationship between the first three moments and cumulants, obtained by extracting 
coefficients from the expansion, is as follows:
\begin{eqnarray}
\kappa_1(T(\rho^\ast))&=&\E[T(\rho^\ast)],\\
\kappa_2(T^2(\rho^\ast))&=&\E[T^2(\rho^\ast)]-(\E[T(\rho^\ast)])^2,\\
\kappa_3(T^3(\rho^\ast))&=&\E[T^3(\rho^\ast)]-3\E[T^2(\rho^\ast)]\E[T(\rho^\ast)]+2(\E[T(\rho^\ast)])^3.
\end{eqnarray}
From the above three relationships (A.8)-(A.10) and (A.5)-(A.7), the first three 
cumulants of $T(\rho^\ast)$ are obtained by
\begin{eqnarray*}
\kappa_1(T(\rho^\ast))&=&\frac{1}{\sqrt{n}}b_1(c)+o\left(n^{-1/2}\right),\\
\kappa_2(T^2(\rho^\ast))&=&\sigma^2(c,a_1,a_2)+o\left(n^{-1/2}\right),\\
\kappa_3(T^3(\rho^\ast))&=&\frac{6}{\sqrt{n}}b_3(c,a_1,a_2,a_3)+o\left(n^{-1/2}\right).
\end{eqnarray*}
Hence, the characteristic function of $T(\rho^\ast)/\sigma(\hat{c},\ah_1,\ah_2)$ can be expressed as 
\begin{eqnarray*}
C(t)&=&
\exp\left(\sum_{j=1}^3\frac{1}{j!}(i t)^j\frac{\kappa_j(T(\rho^\ast))}{\sigma^j(c,a_1,a_2)}\right)+o\left(n^{-1/2}\right)\\
&=&\exp\left(-\frac{t^2}{2}\right)
\left[
1+\frac{1}{\sqrt{n}}\left\{
\frac{b_1(c)(it)}{\sigma(c,a_1,a_2)}
+\frac{b_3(c)(it)^3}{\sigma^3(c,a_1,a_2)}
\right\}
\right]+o\left(n^{-1/2}\right).
\end{eqnarray*}
This result show Proposition 2.2. $\hfill\square$
\subsection*{A.5.~~Proof of Proposition 2.3.}
By using Lemma 2.1, we have
\begin{eqnarray}
&&\Pr\left(\frac{N-p}{p}\frac{T^2}{n}\ge F_{p,N-p}(\alpha)\right)
\to \Phi\left(\frac{\sqrt{n(1-c)}\Delta_{\Sigma^{-1}}^2}{\sqrt{2c}}-z(\alpha)\right),\\
&&\Pr\left(\sqrt{n}\frac{D_n-1}{\sigma_2(\hat{a}_1,\hat{a}_2)}\ge z(\alpha)\right)\to
\Phi\left(\frac{\sqrt{n}\Delta_I^2}{\sqrt{2ca_2}}-z(\alpha)\right),\\
&&\Pr\left(\frac{T(\hat{\rho}^\ast)}{\sigma(c,\ah_1,\ah_2)}\ge z(\alpha)\right)
\to\Phi\left(\sqrt{n}\frac{\sqrt{a_2(1-c)}\De_{\Sigma^{-1}}^2+\De_I^2}{2\sqrt{\{(1-c)^{1/2}a_1a_2^{1/2}+a_2\}c}}-z(\alpha)\right).
\end{eqnarray}
From (A.11) and (A.13), we have
\begin{eqnarray}
& &\Phi\left(\frac{\sqrt{n(1-c)}\Delta_{\Sigma^{-1}}^2}{\sqrt{2c}}-z(\alpha)\right)\leq
\Phi\left(\sqrt{n}\frac{\sqrt{a_2(1-c)}\De_{\Sigma^{-1}}^2+\De_I^2}{2\sqrt{\{(1-c)^{1/2}a_1a_2^{1/2}+a_2\}c}}-z(\alpha)\right)
\notag\\
&\Leftrightarrow&\frac{\sqrt{1-c}\Delta_{\Sigma^{-1}}^2}{\sqrt{2c}}\leq\frac{\sqrt{a_2(1-c)}\De_{\Sigma^{-1}}^2+\De_I^2}{2\sqrt{\{(1-c)^{1/2}a_1a_2^{1/2}+a_2\}c}}
\notag\\
&\Leftrightarrow&\left(\frac{\sqrt{1-c}}{\sqrt{2c}}-\frac{\sqrt{a_2(1-c)}}{2\sqrt{\{(1-c)^{1/2}a_1a_2^{1/2}+a_2\}c}}\right)
\frac{\Delta_{\Sigma^{-1}}^2}{\Delta_{I}^2}
\leq\frac{1}{2\sqrt{\{(1-c)^{1/2}a_1a_2^{1/2}+a_2\}c}}
\notag\\
&\Leftrightarrow&
\frac{\Delta_{\Sigma^{-1}}^2}{\Delta_{I}^2}
\leq\frac{\left\{\sqrt{2}\left(1+a_1\sqrt{(1-c)/a_2}\right)^{1/2}-1\right\}^{-1}}{\sqrt{a_2(1-c)}}.
\end{eqnarray}
Therefore, condition (A.14) are necessary and sufficient conditions for 
the condition that the local asymptotic power of our test is superior to the local asymptotic power of $T^2$-test. 
Similarly, we have that
\begin{eqnarray}
& &\Phi\left(\frac{\sqrt{n}\Delta_I^2}{\sqrt{2ca_2}}-z(\alpha)\right)\leq
\Phi\left(\sqrt{n}\frac{\sqrt{a_2(1-c)}\De_{\Sigma^{-1}}^2+\De_I^2}{2\sqrt{\{(1-c)^{1/2}a_1a_2^{1/2}+a_2\}c}}-z(\alpha)\right)\notag\\
&\Leftrightarrow&\frac{\Delta_I^2}{\sqrt{2ca_2}}\leq\frac{\sqrt{a_2(1-c)}\De_{\Sigma^{-1}}^2+\De_I^2}{2\sqrt{\{(1-c)^{1/2}a_1a_2^{1/2}+a_2\}c}}
\notag\\
&\Leftrightarrow&\left(\frac{1}{\sqrt{2ca_2}}-\frac{1}{2\sqrt{\{(1-c)^{1/2}a_1a_2^{1/2}+a_2\}c}}
\right)
\leq\frac{\sqrt{a_2(1-c)}}{2\sqrt{\{(1-c)^{1/2}a_1a_2^{1/2}+a_2\}c}}\frac{\Delta_{\Sigma^{-1}}^2}{\Delta_{I}^2}\notag\\
&\Leftrightarrow&
\frac{\sqrt{2}\left(1+a_1\sqrt{(1-c)/a_2}\right)^{1/2}-1}{\sqrt{a_2(1-c)}}
\leq\frac{\Delta_{\Sigma^{-1}}^2}{\Delta_{I}^2},
\end{eqnarray}
by using (A.11) and (A.13). Therefore, condition (A.15) are necessary and sufficient conditions for 
the condition that the local asymptotic power of our test is superior to the local asymptotic power of D-test. 
Since $0\leq \sqrt{2}(1+a_1\sqrt{(1-c)/a_2})^{1/2}-1\leq 1$, we have
\begin{eqnarray}
\frac{\sqrt{2}\left(1+a_1\sqrt{(1-c)/a_2}\right)^{1/2}-1}{\sqrt{a_2(1-c)}}\leq 
\frac{\left\{\sqrt{2}\left(1+a_1\sqrt{(1-c)/a_2}\right)^{1/2}-1\right\}^{-1}}{\sqrt{a_2(1-c)}}.
\end{eqnarray}
Combining (A.14)-(A.16), we obtain
\begin{eqnarray*}
&{\rm (i)}&
\frac{\sqrt{2}\left(1+a_1\sqrt{(1-c)/a_2}\right)^{1/2}-1}{\sqrt{a_2(1-c)}}
\leq\frac{\Delta_{\Sigma^{-1}}^2}{\Delta_{I}^2}\leq\frac{\left\{\sqrt{2}\left(1+a_1\sqrt{(1-c)/a_2}\right)^{1/2}-1\right\}^{-1}}{\sqrt{a_2(1-c)}}\\
&\Leftrightarrow&
\lim_{n,p\to\infty}\Pr\left(\frac{T(\hat{\rho}^\ast)}{\sigma(c,\ah_1,\ah_2)}\ge z(\alpha)\right)\\
& &
>\max\left\{\lim_{n,p\to\infty}\Pr\left(\sqrt{n}\frac{D_n-1}{\sigma_2(\hat{a}_1,\hat{a}_2)}\ge z(\alpha)\right),
\lim_{n,p\to\infty}\Pr\left(\frac{N-p}{p}\frac{T^2}{n}\ge F_{p,N-p}(\alpha)\right)\right\},\\
&{\rm (ii)}&
\frac{\sqrt{2}\left(1+a_1\sqrt{(1-c)/a_2}\right)^{1/2}-1}{\sqrt{a_2(1-c)}}
>\frac{\Delta_{\Sigma^{-1}}^2}{\Delta_{I}^2}\\
&\Leftrightarrow&
\lim_{n,p\to\infty}\Pr\left(\sqrt{n}\frac{D_n-1}{\sigma_2(\hat{a}_1,\hat{a}_2)}\ge z(\alpha)\right)>
\lim_{n,p\to\infty}\Pr\left(\frac{T(\hat{\rho}^\ast)}{\sigma(c,\ah_1,\ah_2)}\ge z(\alpha)\right)\\
& &>\lim_{n,p\to\infty}\Pr\left(\frac{N-p}{p}\frac{T^2}{n}\ge F_{p,N-p}(\alpha)\right),\\
\end{eqnarray*}
\begin{eqnarray*}
&{\rm (iii)}&
\frac{\left\{\sqrt{2}\left(1+a_1\sqrt{(1-c)/a_2}\right)^{1/2}-1\right\}^{-1}}{\sqrt{a_2(1-c)}}<\frac{\Delta_{\Sigma^{-1}}^2}{\Delta_{I}^2}\\
&\Leftrightarrow&
\lim_{n,p\to\infty}\Pr\left(\frac{N-p}{p}\frac{T^2}{n}\ge F_{p,N-p}(\alpha)\right)>
\lim_{n,p\to\infty}\Pr\left(\frac{T(\hat{\rho}^\ast)}{\sigma(c,\ah_1,\ah_2)}\ge z(\alpha)\right)\\
& &>\lim_{n,p\to\infty}\Pr\left(\sqrt{n}\frac{D_n-1}{\sigma_2(\hat{a}_1,\hat{a}_2)}\ge z(\alpha)\right).
\end{eqnarray*}
These results prove Proposition 2.3. 
$\hfill\square$

\newpage
\begin{table}[H]
\centering
\begin{center}
${\rm Table}~1.~{\rm ASL~in~the~case~of}~(\eta,p)=(0.2,50)$
\end{center}
\begin{tabular}{ccccccc}
\hline
$\alpha~\backslash~N$ & & $70$ & $110$ & $150$ & $190$ & $230$ \\
\hline
                           & $T^2$ & 0.009 & 0.010 & 0.010 & 0.010 & 0.010\\
$0.01$                     &$D_n$ & 0.011 & 0.010 & 0.010 & 0.010 & 0.010\\
                           &$T(\rho)$ & 0.009 & 0.009 &0.009 & 0.009& 0.009\\
\hline
                           & $T^2$    & 0.046 & 0.049 & 0.050 & 0.051 & 0.050\\
$0.05$                     &$D_n$     & 0.051 & 0.051 & 0.050 & 0.050 & 0.050\\
                           &$T(\rho)$ & 0.043 & 0.046 & 0.047 & 0.048 & 0.048\\
\hline
                           & $T^2$ & 0.095 & 0.099 & 0.100 &0.101 &0.100\\
$0.10$                     &$D_n$ &  0.101 & 0.100 & 0.101 & 0.100 & 0.100\\
                           &$T(\rho)$ & 0.089 & 0.094 & 0.096& 0.097& 0.097\\
\hline\\
\end{tabular}
\centering
\begin{center}
${\rm Table}~2.~{\rm ASL~in~the~case~of}~(\eta,p)=(0.4,50)$
\end{center}
\begin{tabular}{ccccccc}
\hline
$\alpha~\backslash~N$ & & $70$ & $110$ & $150$ & $190$ & $230$ \\
\hline
                           & $T^2$ & 0.009 & 0.010 & 0.010 & 0.009 & 0.011\\
$ 0.01$                    &$D_n$ &0.011 &0.010 & 0.010 & 0.010 & 0.010\\
                           &$T(\rho)$ & 0.009 & 0.008 & 0.009 &0.009 &0.010\\
\hline
                           & $T^2$ & 0.046 & 0.050 & 0.050 & 0.050 & 0.051\\
$0.05$                     &$D_n$ & 0.052 & 0.051 & 0.050 & 0.049 & 0.051\\
                           &$T(\rho)$ &0.044 &0.045 &0.046 &0.046 &0.049\\
\hline
                           & $T^2$ & 0.095 & 0.010 & 0.101 & 0.099 & 0.100\\
$0.10$                     &$D_n$ &0.102 & 0.102 & 0.099 & 0.099 & 0.101\\
                           &$T(\rho)$ & 0.091 & 0.095 &0.095 &0.094 &0.097\\
\hline\\
\end{tabular}
\centering
\begin{center}
${\rm Table}~3.~{\rm ASL~in~the~case~of}~(\eta,p)=(0.6,50)$
\end{center}
\begin{tabular}{ccccccc}
\hline
$\alpha~\backslash~N$ &  & $70$ & $110$ & $150$ & $190$ & $230$ \\
\hline
                           & $T^2$ & 0.009 & 0.010 & 0.009 & 0.011 &0.010\\
$0.01$                     &$D_n$ & 0.010  & 0.010 & 0.010 & 0.009 & 0.010\\
                           &$T(\rho)$ & 0.008 & 0.008 &0.009 &0.010 &0.010\\
\hline
                           & $T^2$ & 0.046 & 0.049 & 0.048 & 0.049 &0.050\\
$0.05$                     &$D_n$  & 0.051 & 0.050 & 0.050 & 0.050 & 0.049\\
                           &$T(\rho)$ &0.043 &0.045 &0.046 &0.047 & 0.047\\
\hline
                           & $T^2$ & 0.095 & 0.100 & 0.099 &0.100 &0.100\\
$0.10$                     &$D_n$  &0.102 & 0.100 & 0.100 & 0.100 & 0.099\\
                           &$T(\rho)$ & 0.089 &0.093 &0.093 &0.096 &0.095\\
\hline\\
\end{tabular}
\end{table}
\begin{table}[H]
\centering
\begin{center}
${\rm Table}~4.~{\rm ASL~in~the~case~of}~(\eta,p)=(0.2,100)$
\end{center}
\begin{tabular}{ccccccc}
\hline
$\alpha~\backslash~N$ & & $120$ & $160$ & $200$ & $240$ & $280$ \\
\hline
                           &$T^2$ & 0.008 & 0.009 & 0.010 & 0.010 & 0.010\\
$0.01$                     &$D_n$ &0.009 & 0.010 & 0.010 & 0.010 & 0.010\\
                           &$T(\rho)$ & 0.008 &0.007 &0.009 &0.008 & 0.009\\
\hline
                           &$T^2$ & 0.046 & 0.048 & 0.049 & 0.050 & 0.050\\
$0.05$                     &$D_n$ & 0.050 & 0.049 & 0.050 & 0.050 & 0.049\\
                           &$T(\rho)$ & 0.041 &0.043 &0.045 &0.046 & 0.047\\
\hline
                           & $T^2$ & 0.094 & 0.096 & 0.098 & 0.101 & 0.099\\
$0.10$                     &$D_n$ &  0.099 &0.099 & 0.100 & 0.099 & 0.099\\
                           &$T(\rho)$ & 0.088 & 0.089 &0.093 &0.097 &0.095\\
\hline\\
\end{tabular}
\centering
\begin{center}
${\rm Table}~5.~{\rm ASL~in~the~case~of}~(\eta,p)=(0.4,100)$
\end{center}
\begin{tabular}{ccccccc}
\hline
$\alpha~\backslash~N$ & & $120$ & $160$ & $200$ & $240$ & $280$ \\
\hline
                           &$T^2$ & 0.008 & 0.010 & 0.010 & 0.010 & 0.010\\
$ 0.01$                    &$D_n$ & 0.009 &0.011 & 0.010 & 0.010 & 0.010\\
                           &$T(\rho)$ & 0.008 & 0.008 & 0.009 &0.009 &0.009\\
\hline
                           & $T^2$ & 0.046 & 0.049 & 0.050 & 0.050 & 0.051\\
$0.05$                     &$D_n$ & 0.051 & 0.050 & 0.051 & 0.050 & 0.051\\
                           &$T(\rho)$ &0.041 &0.044 &0.046 &0.046 &0.048\\
\hline
                           & $T^2$ & 0.094 & 0.098 & 0.100 & 0.099 & 0.102\\
$0.10$                     &$D_n$ &0.101 &0.100 & 0.099 & 0.100 & 0.100\\
                           &$T(\rho)$ & 0.087 &0.092 &0.095 &0.096 &0.097\\
\hline\\
\end{tabular}
\centering
\begin{center}
${\rm Table}~6.~{\rm ASL~in~the~case~of}~(\eta,p)=(0.6,100)$
\end{center}
\begin{tabular}{ccccccc}
\hline
$\alpha~\backslash~N$ &  & $120$ & $160$ & $200$ & $240$ & $280$ \\
\hline
                           & $T^2$ & 0.008 & 0.010 & 0.010 & 0.010 &0.010\\
$0.01$                     &$D_n$ &0.010 & 0.010 & 0.010 & 0.010 & 0.010\\
                           &$T(\rho)$ & 0.008 &0.008 &0.008 &0.008 &0.009\\
\hline
                           & $T^2$    & 0.046 & 0.049 & 0.050 & 0.049 &0.050\\
$0.05$                     &$D_n$     & 0.051 & 0.050 & 0.050 & 0.050 & 0.050\\
                           &$T(\rho)$ &0.041 &0.044 &0.045 &0.046 &0.048\\
\hline
                           & $T^2$ & 0.095 & 0.099 & 0.099 & 0.099 &0.100\\
$0.10$                     &$D_n$ & 0.100 & 0.100 & 0.100 & 0.100 & 0.100\\
                           &$T(\rho)$ & 0.087 &0.091 &0.094 &0.094 & 0.097\\
\hline\\
\end{tabular}
\end{table}
\begin{table}[H]
\centering
\begin{center}
${\rm Table}~7.~{\rm Empirical~powers~with}~(\eta,p)=(0.2,50)$
\end{center}
\begin{tabular}{ccccccc}
\hline
$\alpha~\backslash~N$ & & $70$ & $110$ & $150$ & $190$ & $230$ \\
\hline
                           & $T^2$    & 0.06       & 0.23      & 0.40      & {\bf 0.53} & {\bf 0.65}\\
$0.01$                     &$D_n$     & {\bf 0.17} & 0.25      & 0.32      & 0.39       & 0.46\\
                           &$T(\rho)$ & 0.12       &{\bf 0.29} &{\bf 0.42} & {\bf 0.53} & 0.63\\
\hline
                           & $T^2$    & 0.21       & 0.49       & 0.67      & {\bf 0.78} & {\bf 0.85}\\
$0.05$                     &$D_n$     & 0.30       & 0.41       & 0.50      & 0.57       & 0.64\\
                           &$T(\rho)$ & {\bf 0.33} & {\bf 0.55} &{\bf 0.69} & 0.77       & 0.84\\
\hline
                           & $T^2$    & 0.34      & 0.64      & 0.79      &{\bf 0.87} &{\bf 0.92}\\
$0.10$                     &$D_n$     &{\bf 0.51} & 0.63      & 0.71      & 0.77      & 0.81\\
                           &$T(\rho)$ & 0.48      &{\bf 0.69} &{\bf 0.80} &{\bf 0.87} & 0.91\\
\hline\\
\end{tabular}
\centering
\begin{center}
${\rm Table}~8.~{\rm Empirical~powers~with}~(\eta,p)=(0.4,50)$
\end{center}
\begin{tabular}{ccccccc}
\hline
$\alpha~\backslash~N$ & & $70$ & $110$ & $150$ & $190$ & $230$ \\
\hline
                           & $T^2$    & 0.04      & 0.12      & 0.21       & 0.30      & 0.39\\
$ 0.01$                    &$D_n$     &{\bf 0.15} &{\bf 0.22} & 0.27       & 0.33      & 0.39\\
                           &$T(\rho)$ & 0.09      & 0.21      & {\bf 0.30} &{\bf 0.39} &{\bf 0.47}\\
\hline
                           & $T^2$    & 0.15      & 0.33      & 0.46      & 0.56       & 0.65\\
$0.05$                     &$D_n$     & 0.26      & 0.36      & 0.43      & 0.49       & 0.56\\
                           &$T(\rho)$ &{\bf 0.27} &{\bf 0.44} &{\bf 0.55} &{\bf 0.64}  &{\bf 0.71}\\
\hline
                           & $T^2$    & 0.26      & 0.47       & 0.60      & 0.70      & 0.77\\
$0.10$                     &$D_n$     &{\bf 0.48} & {\bf 0.58} & 0.65      & 0.72      & 0.76\\
                           &$T(\rho)$ & 0.41      & {\bf 0.58} &{\bf 0.68} &{\bf 0.76} &{\bf 0.81}\\
\hline\\
\end{tabular}
\centering
\begin{center}
${\rm Table}~9.~{\rm Empirical~powers~with}~(\eta,p)=(0.6,50)$
\end{center}
\begin{tabular}{ccccccc}
\hline
$\alpha~\backslash~N$ &  & $70$ & $110$ & $150$ & $190$ & $230$ \\
\hline
                           & $T^2$    & 0.04       & 0.13       & 0.22      & 0.31      &{\bf 0.40}\\
$0.01$                      &$D_n$    & {\bf 0.10} & 0.14       & 0.18      & 0.21      & 0.25\\
                           &$T(\rho)$ & 0.08       & {\bf 0.17} &{\bf 0.25} &{\bf 0.33} &{\bf 0.40}\\
\hline
                           & $T^2$    & 0.15      & 0.33      & 0.47      & 0.57      &{\bf 0.66}\\
$0.05$                     &$D_n$     & 0.18      & 0.24      & 0.29      & 0.33      & 0.37\\
                           &$T(\rho)$ &{\bf 0.23} &{\bf 0.39} &{\bf 0.50} &{\bf 0.58} & 0.65\\
\hline
                           & $T^2$    & 0.26      & 0.48      & 0.61      &{\bf 0.70} &{\bf 0.78}\\
$0.10$                     &$D_n$     &{\bf 0.38} & 0.45      & 0.51      & 0.56      & 0.61\\
                           &$T(\rho)$ & 0.37      &{\bf 0.53} &{\bf 0.63} &{\bf 0.70} &0.77\\
\hline\\
\end{tabular}
\end{table}
\begin{table}[H]
\centering
\begin{center}
${\rm Table}~10.~{\rm Empirical~powers~with}~(\eta,p)=(0.2,100)$
\end{center}
\begin{tabular}{ccccccc}
\hline
$\alpha~\backslash~N$ & & $120$ & $160$ & $200$ & $240$ & $280$ \\
\hline
                           &$T^2$ & 0.04 & 0.13 & 0.24 & 0.35 & {\bf 0.45}\\
$0.01$                     &$D_n$ &{\bf 0.14} & 0.18 & 0.22 & 0.26 & 0.30\\
                           &$T(\rho)$ & 0.09 &{\bf 0.19} &{\bf 0.28} &{\bf 0.37} & 0.44\\
\hline
                           &$T^2$ & 0.16 & 0.35 & 0.50 & 0.62 &{\bf 0.71}\\
$0.05$                     &$D_n$ & {\bf 0.29} & 0.35 & 0.40 & 0.45 & 0.50\\
                           &$T(\rho)$ & 0.27 &{\bf 0.44} &{\bf 0.55} &{\bf 0.64} & 0.70\\
\hline
                           & $T^2$ & 0.27 & 0.50 & 0.65 & 0.75 & {\bf 0.82}\\
$0.10$                     &$D_n$ &{\bf 0.48} &0.55 & 0.60 & 0.65 & 0.69\\
                           &$T(\rho)$ & 0.42 & {\bf 0.59} &{\bf 0.69} &{\bf 0.77} &{\bf 0.82}\\
\hline\\
\end{tabular}
\centering
\begin{center}
${\rm Table}~11.~{\rm Empirical~powers~with}~(\eta,p)=(0.4,100)$
\end{center}
\begin{tabular}{ccccccc}
\hline
$\alpha~\backslash~N$ & & $120$ & $160$ & $200$ & $240$ & $280$ \\
\hline
                           &$T^2$     & 0.03      & 0.09      & 0.15      & 0.22      & 0.28\\
$ 0.01$                    &$D_n$     &{\bf 0.14} &{\bf 0.17} & 0.21      & 0.25      & 0.28\\
                           &$T(\rho)$ & 0.07      & 0.15      &{\bf 0.22} &{\bf 0.29} &{\bf 0.35}\\
\hline
                           & $T^2$    & 0.12      & 0.26      & 0.37      & 0.47      & 0.54\\
$0.05$                     &$D_n$     &{\bf 0.27} & 0.32      & 0.37      & 0.42      & 0.46\\
                           &$T(\rho)$ & 0.24      &{\bf 0.37} &{\bf 0.47} &{\bf 0.55} &{\bf 0.61}\\
\hline
                           & $T^2$    & 0.22      & 0.39      & 0.51      & 0.61      & 0.68\\
$0.10$                     &$D_n$     &{\bf 0.46} &{\bf 0.52} & 0.57      & 0.62      & 0.66\\
                           &$T(\rho)$ & 0.38      &{\bf 0.52} &{\bf 0.61} &{\bf 0.68} &{\bf 0.73}\\
\hline\\
\end{tabular}
\centering
\begin{center}
${\rm Table}~12.~{\rm Empirical~powers~with}~(\eta,p)=(0.6,100)$
\end{center}
\begin{tabular}{ccccccc}
\hline
$\alpha~\backslash~N$ &  & $120$ & $160$ & $200$ & $240$ & $280$ \\
\hline
                           & $T^2$    & 0.03      & 0.07            & 0.12      & 0.16      &0.22\\
$0.01$                     &$D_n$     &{\bf 0.11} & {\bf 0.14}      &{\bf 0.16} &{\bf 0.19} & 0.21\\
                           &$T(\rho)$ & 0.06      &0.09             & 0.15      & 0.18      &{\bf 0.27}\\
\hline
                           & $T^2$    & 0.11      & 0.22      & 0.31      & 0.39      &0.46\\
$0.05$                     &$D_n$     &{\bf 0.21} & 0.25      & 0.29      & 0.33      & 0.36\\
                           &$T(\rho)$ &{\bf 0.21} &{\bf 0.32} &{\bf 0.40} &{\bf 0.46} &{\bf 0.52}\\
\hline
                           & $T^2$    & 0.20      & 0.34      & 0.45      & 0.53      & 0.61\\
$0.10$                     &$D_n$     &{\bf 0.39} & 0.44      & 0.49      & 0.53      & 0.56\\
                           &$T(\rho)$ & 0.33      &{\bf 0.46} &{\bf 0.54} &{\bf 0.60} &{\bf 0.66}\\
\hline\\
\end{tabular}
\end{table}
\end{document}